\documentclass[11pt,twoside]{article}
 \usepackage{mathptmx}      
   
\usepackage{amsmath} 
\usepackage{subeqnarray}
\usepackage{amsfonts}
\usepackage{amssymb} 
\usepackage{natbib} 
\usepackage{subfigure} 
\usepackage{setspace}
     
\usepackage[margin=1.25in]{geometry}
\usepackage{graphicx}        
\usepackage[format=hang,width=.75\textwidth]{caption}
\usepackage{soul}

\usepackage{xcolor}
\makeatletter
\newcommand{\globalcolor}[1]{%
  \color{#1}\global\let\default@color\current@color
}
\makeatother
 
\usepackage{color}

 \newtheorem{theorem}{Theorem}

 \newtheorem{remark}{Remark}

\def\P#1{\left(#1\right)}
\newcommand{\B}[1]{\mathbf{#1}}

\newcommand{\p}{\partial}
\newcommand{\R}{\mathbb{R}}
\newcommand{\eps}{\varepsilon}
\newcommand{\vfi}{\varphi}
\renewcommand{\O}{\Omega}

\newcommand{\fer}[1]{(\ref{#1})}
\newcommand{\x}{\mathbf{x}}
\newcommand{\y}{\mathbf{y}}
\newcommand{\z}{\mathbf{z}}
\newcommand{\grad}{\nabla}
\newcommand{\abs}[1]{| #1 |}
\newcommand{\nor}[1]{\| #1 \|}
\newcommand{\qed}{$\Box$}

 \usepackage{fancyhdr}
\begin{document}

\title{On a singular perturbation problem arising in the theory of Evolutionary Distributions
\thanks{First author supported by the Spanish MCINN Project MTM2010-18427}
}



\author{Gonzalo Galiano  \thanks{Dpt. of Mathematics, Universidad de Oviedo,
  c/ Calvo Sotelo, 33007-Oviedo, Spain ({\tt galiano@uniovi.es})}
     \and Yosef Cohen\thanks{Department of Fisheries, Wildlife and Conservation Biology, University of Minnesota,  1980 Folwell Avenue, St. Paul, MN 55108, U.S.A.} }

\date{}

\maketitle

\begin{abstract}
Evolution by Natural Selection is a process by which progeny
inherit some properties from their progenitors with small variation.
These properties are subject to Natural Selection and are called adaptive traits
and carriers of the latter are called phenotypes.
The distribution of the density of phenotypes in a population
is called Evolutionary Distributions (ED).
We analyze mathematical models of the dynamics  of a system of ED.
Such systems are anisotropic in that
diffusion of phenotypes in each population
(equation) remains positive in the directions of their 
own adaptive space
and vanishes in the directions of the other's adaptive space.
We prove that solutions to such systems  exist
in a sense weaker than the usual. 
We develop an algorithm for numerical solutions of such systems.
Finally, we conduct numerical experiments---with a
model in which populations compete---that allow us to 
observe salient attributes of a specific system of ED.

\end{abstract}

\section{Introduction}

The---by now discipline---of Evolution by Natural
Selection \cite[]{Darwin+Wallace+1858},
can be distilled to three principles:
(\emph{i}) inheritance with (\emph{ii}) variation (called mutations) and
(\emph{iii}) natural selection.
Phenotypes are traditionally defined as organisms that display some traits.
Traits may be identified with color of skin, height, speed of running,
length of a molecule of protein, volume of cells and so on.
We use a restricted definition of phenotypes, namely an organism that
displays some trait that is subject to natural selection.
With abuse of the English language,
we call these adaptive traits.
For example, an adaptive trait of Polar bears may be the
thickness of hair in their fur,
where length of hair is presumably inherited.

Now consider the speed of running of prey to be an adaptive trait.
Assume that different values of speed are inherited with some variation
and is subject to natural selection.
Then the speed of running is an adaptive trait.
Three issues of interest emerge:
Firstly, a value of a trait is an attribute of an organism.
Secondly, some variation in these values is not ``recognized'' as such
by natural selection.
In our example of speed of running, 
selection influences survival of those prey that run at, say, $20$\,km/hr
and those that run at $22$\,km/hr equally.
This may be so, even when these two different speeds are inherited.
Ergo, a value of an adaptive trait is in fact a set of values.
And this set identifies a subpopulation made of organisms that
belong to a single phenotype where the population is made of phenotypes
that carry all values of the adaptive trait.
We shall get to the third issue in a moment.

Mathematical models of evolution may be cast via 
differential equations that are either deterministic or stochastic. 
Starting with individual-based models is useful for 
one can construct such models based on first principles. 
In such models, particles may represent single organisms that
 may be classified as phenotypes---they exhibit particular values
of adaptive traits.
Particles may also represent a subpopulation of a phenotypes; 
namely a set of phenotypes that are lumped into a single subpopulation%
---of population of all phenotypes---%
based on the granularity of natural selection.
 Granularity refers to the idea that
in acting upon phenotypes, certain range of values of
adaptive traits are indistinguishable by natural selection.
This range defines
a subpopulation of phenotypes which in fact represents 
the unit of natural selection.

In modeling evolution by natural selection,
one often wishes to switch from individual-based
to PDE-based  approach.
One then begins with an individual-based model, 
and assumes that the number of particles, $N\rightarrow$ $\infty$.
The change of approach%
---from individual-based to PDE-based---%
requires justification. 

In a nutshell, one
begins with a definition of some
property that is represented by a particle.
In our case, such property consists of values
of an adaptive trait.
In the context of evolution by natural selection,
the dynamics of values of such property
are represented by changes in the frequency
of such values in a population. 

Next, one passes from
a space-discrete Lagrangian version of the 
dynamics to space-continuous Eulerian version.
In passing, one is faced with two problems: 

Firstly, showing that the counting measure corresponding to the 
sequence defined 
by the $N$-particles system as $N\to\infty$
converges to a continuous limit measure. 
This convergence must be with respect to average
values of adaptive traits. Under appropriate 
regularity conditions, it can be shown that this limit
satisfies a PDE-based problem.

Secondly, proving that the limiting evolution PDE, 
along with appropriate boundary conditions and initial data,
is a well-posed problem.
Furthermore, that the solution to the problem
inherits the essential properties---such as 
positiveness and boundedness---of the discrete measures, 

The first problem, namely that of showing the convergence to a continuous process,
has been addressed by
\cite{Champagnat2006,Champagnat2008}
for the adaptive evolution of a single populations. 
Similar mathematical problems have been addressed in 
spatial dynamics of population 
\cite[]{Oelschlager1989,Morale2005,Capasso2008,Lachowicz2008}, 
the theory of chemotaxis and phototaxis \cite[]{Stevens2000,Levy2008},  
and flow in porous medium \cite[]{Oelschlager1990,Oelschlager2002}.

The second problem relates to solution of a uniformly parabolic ED.
Two approaches were taken to derive these parabolic equations---one
by \cite{Champagnat2006,Champagnat2008} 
and the other by \cite[]{Cohen+2011,Cohen2013}.
The former applies to a single ED,
the latter to a system of ED.
The approach employed by \cite[]{Cohen+2011,Cohen2013}
stems from the original work by \cite[]{Kimura1965}.

When ED is considered as a system of parabolic ED,
the resulting PDE problem is no longer uniformly parabolic,
and in fact has a peculiar structure.
To understand the lack of uniform parabolicity of the system of PDEs
arising in a multipopulation model,
let us recall the essential assumptions of the model:
(\emph{i}) selection affects the dynamics of coevolving populations
through mortality---that changes the frequency of subpopulations
of phenotypes; and
(\emph{ii}) selection and birth are instantaneous events.
Consequently, selection and birth are independent events.
It follows that mutations, as they should, are random
and mutations in each coevolving population occur in a separate
(from other populations) so-called adaptive space.

For example, populations of prey and predators may evolve each in an 
adaptive space made of a single dimension, say, speed of running. 
Then, through selection, changes in the frequency of phenotypes
on the trait value of prey will, in general, 
induce changes in the frequency of phenotypes along the
trait values of predators and vice versa.
Therefore, population densities of prey and predators will be affected by 
the distribution of phenotypes along trait values of their own 
as well as by those of the other population.

\emph{However}, mutations in each population occur
in the populations' adaptive space. 
Because mutations cause diffusion---along values
of the adaptive trait---population densities depend on trait 
values of both populations but only diffuse with respect to one of them.

In other words, we are faced with a special case of anisotropic diffusion in which
diffusion of each population remains positive in the directions of their own adaptive space
and vanishes in the directions of the other's adaptive space.


%

Recall that the second problem we are faced with is
showing that the limiting evolution PDE, 
along with appropriate boundary conditions and initial data,
is a well-posed problem.
Chipot and coworkers studied related issues
in the framework of the theory of singular perturbation. 
\cite{Chipot2011} proved existence of solutions
for a general class of time-independent boundary value problems. 
They also proved
existence  and some additional properties of solutions for  a
particular case of a nonlocal evolution boundary value problem. 
The problems they studied addressed a single equation;
we address a system of equations.
Yet, we fundamentally rely on their work.

The complexity of the dynamics of evolution implies that exact solutions to the problem in terms of 
elementary functions are not available. Even analytical qualitative properties of solution are hard to obtain. 
Thus, in Section (\ref{sec: Discretization}), we introduce a numerical 
discretization of the problem under the
Finite Element Method framework.
Then, we construct some numerical experiments to observe salient
features of the system of ED that cannot be observed otherwise.
We also use the numerical experiments to draw conclusions about potential
dynamics of evolution under such peculiar systems.
See \cite{Galiano2011,Galiano2011a,Galiano2012}
for related numerical approaches.

\section{Mathematical model}
\label{sec.mathmodel}

Consider a population in which individuals give
birth and die at rates that are determined by values
of an adaptive trait, $\x$, that they
exhibit (a so-called phenotype) 
and by interactions with other phenotypes. 
The population is characterized at any time $t$ by the finite
counting measure

\begin{equation*}
 \nu_t=\sum_{j=1}^{N(t)} \delta_{\x_t^j}
\end{equation*}
where $\delta_\x$ is the Dirac measure at $\x$. 
The measure $\nu_t$
describes the distribution of phenotypes over the trait space
at time $t$.
Here $N(t)$ is the total number of phenotypes alive
at time $t$, and $\x_t^1,\ldots,\x_t^{N(t)}$ denote the 
values of an adaptive trait---which we identify with individuals
(or particles). 

Let $K$ be the initial number of individuals and  define a normalized 
population process by
\begin{equation*}
 X_t^K=\frac{1}{K}\nu_t^K.
\end{equation*}

It can be shown that if the initial condition
$X_0^K=\nu_0^K/K$ converges to a finite deterministic 
measure with density $u_0$ as $K\to\infty$,  
then under suitable assumptions
about the dynamics of 
the discrete process $\nu_t$---see \cite{Champagnat2006} 
for details---the limit of $X_t^K$ converges in law to a deterministic
measure with density $u$ satisfying the following integro-differential equation
\begin{align}
\label{eq.champa}
 \p_t u(t,\x)=
 \big( (1-\mu(\x)) &  b(\x, V*u(t,\x))- d(\x, U*u(t,\x))\big)u(t,\x) \\
 & +   \int_{\R^d}\mu(\y) M(\y,\x-\y)b(\x, V*u(t,\y)) u(t,\y)d\y. \nonumber
\end{align}
Here, $U$ and $V$ denote interaction kernels affecting mortality and reproduction, 
respectively. 
Functions $d$ and $b$ define the death and birth rate, respectively, of 
individuals with trait $\x$, whereas $\mu(\x)$ is the probability that an offspring 
produced by a phenotype with trait-value $\x$ carries 
a value of the adaptive trait that had mutated.
The function $M(\x,\z)$ encapsulates the so-called mutation kernel.
It reflects the fact that trait values of progeny, $\x+\z$,
differ from those of the progenitor with trait value $\x$
(this, in fact, is the definition of mutations).
As noticed by \cite{Champagnat2006}, (\ref{eq.champa}) 
is an extension of Kimura's equation \cite[]{Kimura1965}.

When the values of mutations are small 
and are distributed systematically 
(that is according to some probability density), 
the nonlocal equation (\ref{eq.champa})
 may be approximated by 
\begin{eqnarray}
	 \p_t u(t,\x) -\frac{1}{2} \Delta \big(\mu(\x)\eps^2(\x) u(t,\x)\big) &=&  
	 u(t,\x)\big( (1-\mu(\x))   b(x, V*u(t,\x))- \label{eq.pde}\\
	 & & d(x, U*u(t,\x))\big), \nonumber
\end{eqnarray}
%
where $\eps(\x)\geq\eps_0>0$ is the second moment of $M$.
Observe that, because adaptive traits represent biological attributes,
the set of their values, $\Omega \subset \R^d$, must be bounded.
So (\ref{eq.pde}) is assumed to hold in $(0,T)\times\Omega$, 
for some final time $T>0$. 
Phenotypes with trait-values outside the boundaries cannot survive;
therefore, we prescribe homogeneous Dirichlet boundary conditions
 \begin{equation*}
  u(t,\x)=0 \quad\text{on }(0,T)\times\p\Omega.
 \end{equation*}

When the dynamics include two populations, 
the generalization of (\ref{eq.pde}) leads to  a system of 
partial differential equations for the densities phenotypes 
$u_1$ and $u_2$ with respect to their
corresponding populations. 
 It is here that some subtleties arise, and we introduce them next.
But first, some notation.

The open and bounded trait space of each population ($i=1,2$) 
is denoted by $\Omega_i$, with $\Omega_i\subset\R^{d_i}$. 
An element of $\Omega_i$ is denoted by $\x_i$ and  is given by its components 
$\x_i=(x_i^1,\ldots,x_i^{d_i})$. 
We then write $\O=\O_1\times\O_2$ and $\x=(\x_1,\x_2)$.
 We use the Laplacian operator restricted to $\Omega_i$ and write 
\[
\Delta_{\x_i} u(t,\x) := \sum_{j=1}^{d_i}\p_{x_j}^2 u(t,\x) .
\]
We now state the problem for two populations thus. 

For $i,j=1,2$ and $i\neq j$, find $u_i:(0,T)\times \O \to \R_+$ such that 
\begin{align}
& \p_t u_i  -  \Delta_{\x_i} \big( c_i  u_i\big) +F_i(\cdot,\cdot,u_1,u_2)=0
  & &\text{in } (0,T)\times\O ,&  \label{eq: the problem} \\
&u_i=0 & & \text{on }(0,T)\times \p\O_i \times \O_j,& \label{eq.bc} \\
&u_i(0,\cdot)=u_{i0} & & \text{in }\O. & \label{eq.id}
\end{align}
The magnitude and rate of mutations are defined thus:
\begin{equation}
\label{eq: rate of mutation}
 c_i(\x_i)=\frac{1}{2}\mu_i(\x_i) \eps_i^2(\x_i)
\end{equation}
and the interactions within and among populations---including
natural selection---thus
\begin{align}
\label{eq: F}
F_i(t,\x,u_1(t,\x),u_2(t,\x)) = & -\big((1-\mu_i(\x_i)) u_i(t,\x) b_i(\x, V*u_i(t,\x)) \\
  & - u_i(t,\x) d_i(\x, U*u_1(t,\x),U*u_2(t,\x)) \big).\nonumber
\end{align}
Equations (\ref{eq: the problem}) to (\ref{eq: F}) 
specify the problem, which we dub $\B P$. 

\begin{remark}
From a direct differential approach, the following
model was derived in \cite{Cohen+2011}. For $i=1,\ldots,N$ (number of populations), 
find $U_i:(0,T)\times \O \to \R_+$ such that
\begin{equation}
\label{eq: ED}
\partial_t U_{i}\left(\mathbf{x},t\right) - \Delta_{\x_i} \Big(\tilde \mu_i(\x_i) B_i(U_i(t,\x)) \Big)=
  B_i(U_i(t,\x)) - M_i(U_1(t,\x),U_2(t,\x)), 
\end{equation}
for some monotone functions $B_i$, and with auxiliary data similar to \fer{eq.bc}-\fer{eq.id}. Observe that for a suitable choice of $\tilde\mu_i$, $B_i$ and $M_i$, Equation~\fer{eq: the problem} may be recasted in the form \fer{eq: ED}.   The solution of (\ref{eq: ED}) was dubbed Evolutionary Distribution
\cite[]{Cohen+2005}.
\end{remark}

\section{The main result}

In this section we study the well-posedness of a variant of $\B P$,
which we christen $\B P1$,
namely:
\begin{align}
 & \p_t u_i - c_i\Delta_{\x_i} u_i+F_i(\cdot,\cdot,u_1,u_2)=0 & & 
 \quad\text{in }Q_T=(0,T)\times\O, & \label{eq.spde}\\
 & u_i= 0 & & \quad\text{on }(0,T)\times\p\O, & \label{eq.x}\\
 & u_i(0,\cdot)= u_{i0} & & \quad\text{in }\O, & \label{eq.sid}
\end{align}
for $i=1,2$, where we have used the notation introduced in 
Section \ref{sec.mathmodel}. 

We rely on the following hypotheses about the data, 
which we shall refer to as \textbf{H}: 

For $i=1,2$,
\begin{enumerate}
\item  $T\in\R_+$ is arbitrarily fixed, and $\Omega_i\subset\R^{d_i}$ is a bounded set with Lipschitz continuous boundary $\partial \Omega_i$.
\item The diffusion coefficients, $c_i$, are positive constants.
\item The birth-selection terms $F_i:Q_T \times \R^2_+\to\R$ satisfy:
\begin{enumerate}
\item $F_i(\cdot,\cdot,s_1,s_2) \in L^\infty (Q_T)$ for all $(s_1,s_2)\in\R^2_+$. 

\item $F_i(t,\x,\cdot,\cdot) $ is Lipschitz continuous in bounded sets of $\R^2_+$ for a.e.  $(t,\x)\in Q_T$, 
with $F_1(\cdot,\cdot,0,\cdot)=F_2(\cdot,\cdot,\cdot,0)=0$.
 \end{enumerate}
\item The initial data are non-negative and satisfy 
$u_{i0}\in L^\infty(\Omega)\cap H^1_{i0}(\O)$, where 
\begin{equation}
\label{def.space}
 H^1_{i0}(\O)=\left\{\vfi\in H^1(\O) : \vfi=0\text{ on }\O_i \right\}. 
\end{equation}
\end{enumerate}

We construct a solution of problem $\B P1$, namely $(u_1,u_2)$, 
as a limit of a sequence of solutions of approximated singular perturbed problems, 
$(u_1^{(\eps)},u_2^{(\eps)})$. 
In doing so, we face two main difficulties. 

Firstly, the lack of uniform estimates
for $\grad_{\x_j}u_i^{(\eps)}$ for $i\neq j$, 
which prevent us from getting a global 
(in space) formulation of a weak solution of  $\B P1$. 
As pointed out by \cite{Chipot2009}, 
the regularity of $u_i\in L^2(0,T;H^1_0(\O))$ 
can not be guaranteed and therefore, 
the  trace of $u_i$ on the boundary $\O_i\times\p\O_j$
is in general not well defined. 
For examples, see \cite{Chipot2009}, where
a solution  of a
time-independent problem is shown to satisfy 
$u_i\in H^1_0(\O_i)$ concomitant with $u_i \notin H^1_0(\O)$. 
This difficulty also motivates
the introduction of the space $H_{i0}^1(\O)$ in (\ref{def.space}).
Thus, the notion of weak solution is stated in \emph{slices} of the domain
(see Theorem~\ref{th.existence}). 

Secondly, and, again, due to the lack of uniform estimates for the gradients 
in the whole domain $\O$, strong convergence of the sequence 
$(u_1^{(\eps)},u_2^{(\eps)})$ in $L^2(Q_T)\times L^2(Q_T)$
can not be achieved by the usual compactness argument. 
Yet, this type of convergence is needed 
to identify the limits of $F_i(\cdot,\cdot,u_1^{(\eps)},u_2^{(\eps)})$. Thus, we resort
to a monotonicity property implied by the Lipschitz continuity (H)$_{3(b)}$ which allows us to perform this identification. 

\begin{remark}
 
%
\noindent (i) 
With minor modifications, Theorem~\ref{th.existence} may be proved for $\x_i$-dependent diffusion coefficients $c_i$ under regularity assumptions and the conditions $c_i>c_0$, for some positive constant $c_0$, and  
$\Delta_{\x_i} c_i \leq 0$ in $\O_i$. Indeed, the monotonicity of the elliptic operator in $H^1_0(\O_i)$ holds under these conditions:
\begin{align*}
 -\int_{\O_i} \big(\Delta_{\x_i} (c_iu_i) -\Delta_{\x_i} (c_i \tilde u_i)\big)(u_i-\tilde u_i)& =
 \int_{\O_i} c_i \abs{\grad_{\x_i} (u_i-\tilde u_i)}^2 - \frac{1}{2}\int_{\O_i} \abs{u_i-\tilde u_i}^2 \Delta_{\x_i}c_i   \\
 & \geq c_0 \int_{\O_i}  \abs{\grad_{\x_i} (u_i-\tilde u_i)}^2.
\end{align*}

\noindent (ii) Although to prove our main result on existence of solutions 
we only need to consider a general
Lipschitz domain $\O$, in the context of our application this domain is 
in fact a regular polyhedra due to the choice of independent traits.
\end{remark}

In the following theorem we assert the existence of solutions of 
problem (\ref{eq.spde})-(\ref{eq.sid})
in a sense weaker than the usual. 
Indeed, although the components of the solution are globally bounded 
in $Q_T$, the existence of their weak derivatives may be justified 
only in slices and so the notion of weak
solution.

\begin{theorem}
 \label{th.existence}
 Assume   {\bf H}. Then problem (\ref{eq.spde})-(\ref{eq.sid}) 
 has a weak solution $(u_1,u_2)$ satisfying, 
 for $i=1,2$,  $u_i\geq 0$ in $Q_T$, 
\begin{equation}
\label{regularity.th}
u_i \in L^\infty(Q_T)\cap L^2(0,T;H^1_0(\Omega_i)) \cap H^1(0,T;L^2(\Omega_i)), 
\end{equation}
and, for all $\vfi\in L^{2}(0,T;H^1_0(\Omega_i))$ and a.e. $\x_j\in\O_j$ ($j\neq i$)
\begin{align}\label{weak}
\int_0^T \int_{\O_i} \p_t u_i \, \vfi  
+ c_i \int_0^T \int_{\O_i} \grad_{\x_i} u_i \cdot\grad_{\x_i}\vfi  
+\int_0^T \int_{\O_i} F_i(u_1,u_2)\,\vfi=0 ,
\end{align}
with the initial data being satisfied in the sense that
\begin{align*}
u_1(0,\cdot,\x_2)=u_{10}\quad\text{ in $\O_1$ for a.e. }\x_2\in\O_2, \\  
 u_2(0,\x_1,\cdot)=u_{20}\quad\text{ in $\O_2$ for a.e. }\x_1\in\O_1.
\end{align*}
\qed
\end{theorem}

There is no general theory of existence of solutions for problem  
(\ref{eq.spde})-(\ref{eq.sid}).
Therefore, we introduce the following singular perturbation approximation to the problem. 
We follow an approach similar to that taken by
\cite{Chipot2011} and \cite{Chipot2011a}.
 
\noindent\emph{Proof of Theorem~\ref{th.existence}.} 
We divide the proof in two steps, according to the monotone behavior of functions $F_i$.

\noindent\textbf{Step 1.} Let us suppose that $F_i$ satisfy the following property:
for all $(s_1,s_2)$ and $(\sigma_1,\sigma_2)$ in $\R^2_+$,
\begin{align}
 \label{prop.mon}
  \big(F_1(\cdot,\cdot,s_1,s_2)  - &  
  F_1(\cdot,\cdot,\sigma_1,\sigma_2)\big) (s_1-\sigma_1)+  \\ 
 & \big(F_2(\cdot,\cdot,s_1,s_2)-F_2(\cdot,\cdot,\sigma_1,
 \sigma_2)\big)(s_2-\sigma_2) \geq 0 \nonumber
 \end{align}
 a.e. in $Q_T$. Then we proceed as follows.

For $i,j=1,2$, $i\neq j$, and some small $\eps>0$, 
we wish to find $(u_1^{(\eps)},u_2^{(\eps)})$ such that 
\begin{align}
\int_0^T<\p_t u_i,\vfi>  
+ c_i \int_{Q_T} \grad_{\x_i} u_i \cdot\grad_{\x_i}\vfi & +\eps \int_{Q_T} 
\grad_{\x_j} u_i \cdot\grad_{\x_j}\vfi  \label{weakeps} \\ 
& + \int_{Q_T} F_i(\cdot,\cdot,u_1,u_2)\,\vfi=0 , \nonumber
\end{align}
for all $\vfi\in L^{2}(0,T;H^1_0(\Omega))$, and satisfying 
the initial data 
$u_{i}^{(\eps)}(0,\cdot)=u_{i0}^{(\delta)}$, where 
$u_{i0}^{(\delta)} \in H_0^1(\O)$ and $u_{i0}^{(\delta)}\to u_{i0}$ 
strongly in $L^2(\O_j;H_{0}^1(\O_i))$ as $\delta\to 0$. 
Up to the last part of the
proof, we will consider the parameter $\delta$ fixed,
and for clarity, we will not reference it.

The existence of non-negative solutions of the semilinear evolution problem 
(\ref{weakeps}) have been established for a large range of functions $F_i$, 
including those satisfying (H)$_3$ and (\ref{prop.mon}) 
\cite[see][Theorem 9 in Chapter 7]{Friedman1964}. 
 The following regularity holds \cite[][Theorem 5 in Chapter 7]{Evans1998}:
 \begin{equation}
  \label{regularity}
 \begin{array}{l}
  u_i^{(\eps)}\in L^\infty(Q_T)\cap C([0,T];L^2(\O))\cap L^2(0,T;H^2(\Omega))\cap L^\infty(0,T;H_0^1(\Omega)) ,\\
  \p_t u_i^{(\eps)}\in L^2(Q_T). 
 \end{array}
 \end{equation}
Recall that, due to property (H)$_{3\text{(b)}}$, the $L^\infty(Q_T)$ bound implied by (\ref{regularity}) is uniform with respect to $\eps$ 
\cite[][Theorem 7 in Chapter 7]{Friedman1964}.

We use $\vfi=u_i^{(\eps)}$ as a test 
function in (\ref{weakeps}) to get, for $t\in(0,T]$,
\begin{align*}
\frac{1}{2}\int_\Omega \abs{u_i^{(\eps)}(t,\cdot)}^2  
 +c_i \int_{Q_t} \abs{\grad_{\x_i} u_i^{(\eps)}}^2 
& +\eps \int_{Q_t} \abs{\grad_{\x_j} u_i^{(\eps)}}^2 \\
& + \int_{Q_t} F_i(u_1^{(\eps)},u_2^{(\eps)})u_i^{(\eps)}=\frac{1}{2}\int_\Omega \abs{u_{i0}^{(\delta)}}^2 .
\end{align*}
Summing for $i=1,2$ and using (\ref{prop.mon})'s property of monotonicity, 
we attain uniform estimates for
\begin{align}\label{est.uniform}
 \nor{u_i^{(\eps)}}_{L^\infty(0,T;L^2(\Omega))},\quad \nor{\grad_{\x_i}u_i^{(\eps)}}_{L^2(Q_T)}  
 ,\quad \sqrt{\eps}\nor{\grad_{\x_j}u_i^{(\eps)}}_{L^2(Q_T)},
\end{align}
implying the existence of some $u_i\in L^2(Q_T)$ such that, for a subsequence (not relabeled)
we have that
\begin{equation}
 \label{eq.convergence}
 u_i^{(\eps)} \to u_i ,\quad
 \grad_{\x_i} u_i^{(\eps)} \to \grad_{\x_i} u_i ,\quad 
 \eps \grad_{\x_j} u_i^{(\eps)} \to 0 
\end{equation}
weakly in $L^2(Q_T)$, and
\begin{equation}\label{conv.linf}
 u_i^{(\eps)} \to u_i \quad\text{weakly $*$ in }  L^\infty(Q_T).
\end{equation}

In addition, by smoothing-approximation procedure on $\p_t u_i^{(\eps)}$, we 
find (see \cite[][Theorem 5 in Chapter 7]{Evans1998})
\begin{align}\label{est.time}
\nor{\p_t u_i^{(\eps)}}_{L^2(Q_T)} \leq 
C \big( \nor{F_i(\cdot,\cdot,u_1^{(\eps)},u_2^{(\eps)}}_{L^2(Q_T)} )
+\nor{u_{i0}^{(\delta)}}_{H^1(\Omega))}\big),
 \end{align}
for some constant $C>0$ independent of $\eps$. In view of (\ref{est.uniform}) and
the continuity of $F_i$ stated in (H)$_{3(\text{b})}$, we also find that  
 $\nor{\p_t u_i^{(\eps)}}_{L^2(Q_T)}$ is uniformly bounded, implying 
\begin{equation}
\label{eq.convergence2}
  \p_t u_i^{(\eps)} \to \p_t u_i  \quad\text{weakly in }  L^2(Q_T).
\end{equation}
Observe that (\ref{eq.convergence}) and (\ref{eq.convergence2})
permit us to pass to the limit $\eps\to0$ in 
the linear terms  of (\ref{weakeps}), 
with the identification $u_i=\lim_{\eps\to0}u_i^{(\eps)}$. 
However, 
this is not the case for the nonlinear terms involving $F_i(\cdot,\cdot,u_1^{(\eps)},u_1^{(\eps)})$, for which a stronger sense 
of convergence is needed to identify their limits. In any case, 
the $L^\infty(Q_T)$ uniform bounds of $u_i^{(\eps)}$ imply, see \cite{Evans1990},
the existence of $\tilde F_{i}\in L^\infty(Q_T)$  such that
\begin{align}\label{conv.linf2}
 F_i(\cdot,\cdot,u_1^{(\eps)},u_2^{(\eps)}) \to \tilde F_{i}, 
 \end{align}
weakly $*$ in $L^\infty(Q_T)$. Taking the limit $\eps\to 0$ in (\ref{weakeps}) we find that $u_i$ satisfies
\begin{align}
\int_{Q_T} \p_t u_i \,\vfi  
+ c_i \int_{Q_T} \grad_{\x_i} u_i \cdot\grad_{\x_i}\vfi  +
 \int_{Q_T} \tilde F_i\,\vfi=0 , \label{conv.inter}
\end{align}
for all $\vfi\in L^{2}(0,T;H^1_0(\Omega))$. 

To deduce the convergence $u_i^{(\eps)}\to u_i$ in a stronger sense, 
we apply the monotonicity that is afforded us by 
assumption (\ref{prop.mon}). 
Let 
\begin{align*}
I_i^{(\eps)}=& \int_{Q_t} \Big( (u_i^{(\eps)}-u_i) \p_t (u_i^{(\eps)}-u_i)+
c_i \grad_{\x_i}(u_i^{(\eps)}-u_i)\cdot\grad_{\x_i}(u_i^{(\eps)}-u_i)\\
& + \eps  \grad_{\x_j}u_i^{(\eps)}\cdot\grad_{\x_j}u_i^{(\eps)}+
(u_i^{(\eps)}-u_i)\big(F_i(\cdot,\cdot,u_1^{(\eps)},u_2^{(\eps)})
-F_i(\cdot,\cdot,u_1,u_2)\big) \Big).
\end{align*}
On the one hand, using (H)$_{\text{3(c)}}$ we find 
\begin{align}
\label{est.inf2}
 I_1^{(\eps)}+I_2^{(\eps)}\geq \sum_{i=1}^2 \Big( \int_\O \big(\abs{u_i^{(\eps)}(t,\cdot)-u_i(t,\cdot)}^2
 + c_i \int_{Q_t} \abs{\grad_{\x_i}(u_i^{(\eps)}-u_i)}^2 
 + \eps \int_{Q_t} \abs{\grad_{\x_j}u_i^{(\eps)}}^2\Big).
\end{align}
On the other hand, expanding $I_i^{(\eps)}$ and using the problem
satisfied by $(u_1^{(\eps)},u_2^{(\eps)})$, i.e. (\ref{weakeps}), we find
\begin{align*}
I_i^{(\eps)}=\int_{Q_t} \Big( -u_i^{(\eps)} \p_t u_i-u_i \p_t u_i^{(\eps)}+
u_i\p_t u_i-
2c_i \grad_{\x_i}u_i^{(\eps)}\cdot \grad_{\x_i} u_i+ c_i \grad_{\x_i}u_i\cdot \grad_{\x_i} u_i\\
-u_i^{(\eps)} F_i(\cdot,\cdot,u_1,u_2)
-u_i F_i(\cdot,\cdot,u_1^{(\eps)},u_2^{(\eps)})
+u_i F_i(\cdot,\cdot,u_1,u_2)\Big).
\end{align*}
Taking $\eps\to 0$ and using the convergence properties (\ref{eq.convergence}) and (\ref{eq.convergence2}), and the problem satisfied by $(u_1,u_2)$, i.e.  (\ref{conv.inter}), we obtain $
I_i^{(\eps)}\to 0$.
Therefore, we deduce from (\ref{est.inf2}) the convergences 
\begin{equation}
 u_i^{(\eps)}\to u_i, \quad\grad_{\x_i}u_i^{(\eps)}\to \grad_{\x_i}u_i
\quad\text{and}\quad \eps\grad_{\x_j}u_i^{(\eps)}\to 0,
 \end{equation}
strongly in $L^2(Q_T)$. We may now go back to (\ref{conv.linf2}) and identify  $\tilde F_i$ as the strong limit in $L^2(Q_T)$ of 
$F_i(\cdot,\cdot,u_1^{(\eps)},u_2^{(\eps)})$, as $\eps\to 0$. Hence, (\ref{conv.inter})
becomes 
\begin{align*}
\int_{Q_T}\p_t u_i \,\vfi  
+ c_i \int_{Q_T} \grad_{\x_i} u_i \cdot\grad_{\x_i}\vfi  +
 \int_{Q_T} F_i(\cdot,\cdot,u_1,u_2)\,\vfi=0 , 
\end{align*}
for all $\vfi\in L^{2}(0,T;H^1_0(\Omega))$. Observe that due to (\ref{eq.convergence2}) we also have $u_{i}(0,\cdot)=u_{i0}^{(\delta)}$ a.e. in $\O$.

Finally, we justify the passing to the limit $\delta\to 0$.
Consider test functions of the form $\vfi(t,\x)=\vfi_i(t,\x_i)\vfi_j(\x_j)$, with $\vfi_i\in L^2(0,T;H_{0}^1(\O_i))$ and $\vfi_j\in H^1_0(\O_j)$, for $i,j=1,2$, $i\neq j$. By density, we deduce that a.e. in $\O_j$,
\begin{align*}
\int_0^T\int_{\O_i}\p_t u_i \,\vfi_i  
+ c_i \int_0^T\int_{\O_i} \grad_{\x_i} u_i \cdot\grad_{\x_i}\vfi_i  +
 \int_0^T\int_{\O_i} F_i(\cdot,\cdot,u_1,u_2)\,\vfi_i=0  
\end{align*}
for all $\vfi_i\in L^2(0,T;H_{0}^1(\O_i))$. 
We now may use $\vfi_i=u_i(\cdot,\cdot,\x_j)$ and (a regularization of) 
$\vfi_i=\p_t u_i(\cdot,\cdot,\x_j)$, for a.e. $\x_j\in\O_j$, to obtain similar estimates to those for the problem
 stated in the whole $\O$, e.g. (\ref{est.uniform}) 
and (\ref{est.time}). 

These estimates may be obtained independently of $\delta$ due to 
the regularity $u_{i0}\in H^1_{i0}(\O)$ given in assumption (H)$_4$,
and the convergence $u_{i0}^{(\delta)}\to u_{i0}$ 
strongly in $L^2(\O_j;H_{0}^1(\O_i))$. Since the monotonicity properties 
hold as well, the passing to the limit and its identification follows. 

\noindent\textbf{Step 2.} If $F_i$ do not satisfy property (\ref{prop.mon}), then we introduce the following auxiliary
problem.
Let $\lambda>0$ be a constant to be chosen. Under the change  of unknowns 
 $u_i(t,x)=\text{e}^{\lambda t}U_i(t,x)$  problem (\ref{eq.spde})-(\ref{eq.sid}) is formally equivalent to
\begin{align}
 & \p_t U_i - \Delta_{\x_i} U_i+
 \hat F_i(\cdot,\cdot,U_1,U_2)=0 & & \quad\text{in }Q_T, & \label{eq:aux}\\
 & U_i= 0 & & \quad\text{on }(0,T)\times\p\Omega, & \label{bc:aux} \\
 & U_i(0,\cdot)= u_{i0} & & \quad\text{in }\Omega, & \label{id:aux}
\end{align}
with $\hat F_i(t,\x,s_1,s_2)= \lambda s_i +
\text{e}^{-\lambda t}F_i(t,\x,\text{e}^{\lambda t}s_1,\text{e}^{\lambda t}s_2)$. 
Let $L$ be an upper bound for the Lipschitz constants of $F_i$. 
Using the Lipschitz continuity (H)$_2$, we obtain 
\begin{align*}
  &\big(\hat F_1(s_1,s_2)  -   \hat F_1(\sigma_1,\sigma_2)\big) (s_1-\sigma_1)+  
  \big(\hat F_2(s_1,s_2)-\hat F_2(\sigma_1,\sigma_2)\big)(s_2-\sigma_2)\geq\\ 
  &  
 \lambda \big( (s_1-\sigma_1)^2+(s_2-\sigma_2)^2 \big)-
 L \big( (s_1-\sigma_1)^2+(s_2-\sigma_2)^2 + 
 2\abs{s_1-\sigma_1}\abs{s_2-\sigma_2} \big), 
 \end{align*}
which is non-negative for $\lambda$ large enough. Therefore,  $\hat F_i$ satisfy the monotonicity
condition (\ref{prop.mon}).

Step 1 of this proof 
ensures the existence  of a solution $(U_1,U_2)$ to problem (\ref{eq:aux})-(\ref{id:aux}) in the weak sense 
of (\ref{weak}), and with the 
regularity (\ref{regularity.th}). Then, it is straightforward 
to check that $u_i(t,x)=\text{e}^{\lambda t}U_i(t,x)$  solves problem (\ref{eq.spde})-(\ref{eq.sid}) in an 
identical weak sense.
\qed

\section{Applications and numerical examples}
\label{sec: Discretization}

\subsection{Discretization scheme}

We follow the ideas of Theorem~\ref{th.existence} to construct numerical approximations to 
problem \fer{eq.spde}-\fer{eq.sid}. We first consider the $\eps-$regularized version given by 
problem \fer{weakeps}, a semilinear parabolic problem for which standard 
discretization methods may be used, and then produce simulations of the 
solutions $(u_1^{(\eps)},u_2^{(\eps)})$ for a collection of decreasing 
values of $\eps$ to get an idea of the behaviour of solutions to the original 
unperturbed problem. 

For the numerical approximation, we use finite elements in space and 
backward finite differences in time. 
More concretely, we consider  a quasi-uniform mesh of a rectangle
 $\Omega\subset\R^2$, $\{\mathcal{T}_h\} _h$, where  $h$ 
 represents triangle diameter, and 
introduce the finite element space of piecewise $\mathbb{P}_1$-elements:
$$
S^h = \{ \chi\in \mathcal{C}(\overline{\Omega} ) ; \, 
\chi |_K\in\mathbb{P}_1\,\text{ for all } K\in\mathcal{T}_h \} .
$$
The Lagrange interpolation operator, denoted by 
$\Pi ^h : \mathcal{C}(\overline{\Omega} ) \to S^h$, from 
which we define the discrete semi-inner product on $\mathcal{C}(\overline{\Omega} ) $ 
 is given by 
$$
(\eta_1,\eta_2)^h= \int_{\Omega} \Pi^h(\eta_1,\eta_2) .
$$
For the time discretization, we  use a uniform partition of $[0,T]$ 
of time step $\tau$.
For $t=t_0=0$, set $u_{i0}^{(\eps)}=u_{i0}$.
Then, for $n\geq 1$ find $u_{i,n}^{(\eps)}$ such that for 
$i,j=1,2$ and $i\neq j$,  
\begin{equation}\label{eq.fem}
\begin{array}{l}
\frac{1}{\tau}\big( u_{i,n}^{(\eps)}-u_{i,n-1}^{(\eps)} , \chi )^h
+ \Big( \big( c_i\p_{x_i}u_{i,n}^{(\eps)}, \eps\p_{x_j} u_{i,n}^{(\eps)}\big) ,
\grad\chi \Big) +\big( F_i(u_{1,n-1}^{(\eps)},u_{2,n-1}^{(\eps)}   ),\chi \big)^h =0,
\end{array}
\end{equation}
for every $ \chi\in S^h $.
Since \fer{eq.fem} is a nonlinear algebraic problem, we use a fixed point argument to 
approximate its solution,  $(u_{1,n}^{(\eps)},u_{2,n}^{(\eps)}   )$, 
at each time slice $t=t_n$, from the previous
approximation $(u_{1,n-1}^{(\eps)},u_{2,n-1}^{(\eps)}   )$.  
Let $u_{i,n}^{\eps,0}=u_{i,n-1}^{(\eps)}$. 
Then, for $k\geq 1$ the problem is to find $u_{i,n}^{\eps,k}$ such that for 
 all $\chi \in S^h$ 
\begin{equation*}
\begin{array}{l}
 \displaystyle\frac{1}{\tau}\big( u_{i,n}^{\eps,k}-u_{i,n-1}^{(\eps)} , \chi )^h
 +\Big( \big( c_i\p_{x_i}u_{i,n}^{\eps,k}, \eps\p_{x_j} u_{i,n}^{\eps,k}\big) ,
 \grad\chi \Big) +\big( F_i(u_{1,n}^{\eps,k-1},u_{2,n}^{\eps,k-1}   ),\chi \big)^h =0.
\end{array}
\end{equation*}
We use a stopping criteria based on the $L^2$ relative error between two iterations, i.e.
 \begin{equation*}
 \Big(\sum _{i=1,2} \nor{u_{i,n}^{\eps,k}-u_{i,n}^{\eps,k-1}}_{L^2}^2\Big)^{1/2}
 \Big(\sum _{i=1,2} \nor{u_{i,n}^{\eps,k-1}}_{L^2}^2\Big)^{-1/2}
 <\text{tol},
 \end{equation*}
for values of $\text{tol}$ chosen empirically, and set $u_{i,n}^{(\eps)}=u_{i,n}^{\eps,k}$.
In some of the experiments we integrate in time until a numerical stationary 
solution is achieved. This is determined again by an $L^2$ relative error measure:
 \begin{equation*}
 \Big(\sum _{i=1,2} \nor{u_{i,n}^{(\eps)}-
 u_{i,n-1}^{(\eps)}}_{L^2}^2\Big)^{1/2}  
 \Big(\sum _{i=1,2} \nor{u_{i,n-1}^{(\eps)}}_{L^2}^2\Big)^{-1/2}<\text{tol}_S,
 \end{equation*}
where $\text{tol}_S$ is chosen empirically. 

Unless otherwise stated, we use the parameter values of Table~\ref{table.parameter} 
in all the experiments. The spatial domain is the square $\O=(0,1)\times(0,1)$.
\begin{table}
\centering
\caption{Parameter values for the experiments} 
\begin{tabular}{|l|c|c|c|}
\hline\hline 
Parameter & Symbol & Experiment 1   \\ \hline\hline
Number of spatial nodes & $N$ & $900$ \\
Time step & $\tau$ & $1.e-03$  \\
Initial densities & $(u_{10},u_{20})$ & $(0.5,0.5)$ \\
Fixed point tolerance & tol & $1.e-03$ \\
Stationary state tolerance& $\text{tol}_S$ & $1.e-05$  \\
Growth$^*$ &$(\alpha_1,\alpha_2)$ & $(5,4)$  \\
Selection$^*$ &$\beta_{ij}$ & $\left(\begin{array}{cc} 3& 2 \\ 2 & 2\end{array}\right)$  \\
\hline\hline
\end{tabular}
\label{table.parameter}
\\
\small (*)  $F_i(u_1,u_2)=-\big(\alpha_i u_i + u_i(\beta_{i1}u_1 +  \beta_{i2}u_2)\big)$
\end{table}

\subsection{Numerical examples}

Here we show the outcome of competition between
two distributions of phenotypes, $u_1$ and $u_2$, for various values of diffusion 
coefficients $c_1$ and $c_2$ and for interactions within and among populations given by 
\begin{equation*}
	F_i\P{u_1,u_2} = -\Big(\alpha_i u_i + u_i \P{\beta_{i1} u_1 + \beta_{i2}u_2}\Big),
\end{equation*}
see Table~\ref{table.parameter} for parameter values. 

In Experiments~1 and 2 we explore the behavior of solutions for the cases $c_1 \gg c_2$ and $c_1=c_2$, respectively,
directly implementing the discrete version of problem P1, i.e. with the perturbation parameter $\eps=0$. 
Although Theorem~\ref{th.existence} ensures existence of solutions defined in a weaker sense than the usual, 
we observe good regularity properties in their numerical approximations. 

We explore the effect of the magnitude of $\varepsilon$
on the singular perturbation approximation, i.e. to the outcome of \fer{weakeps}.
Thus, in Experiment~3 we produce numerical approximations of \fer{weakeps}
for several decreasing values of $\eps$.

\subsubsection{Experiment 1: $c_1 = 0.1, c_2 = 0.01, \varepsilon = 0$}
 
\begin{figure}[t]
\centering
  {\includegraphics[width=6.cm,height=6cm]{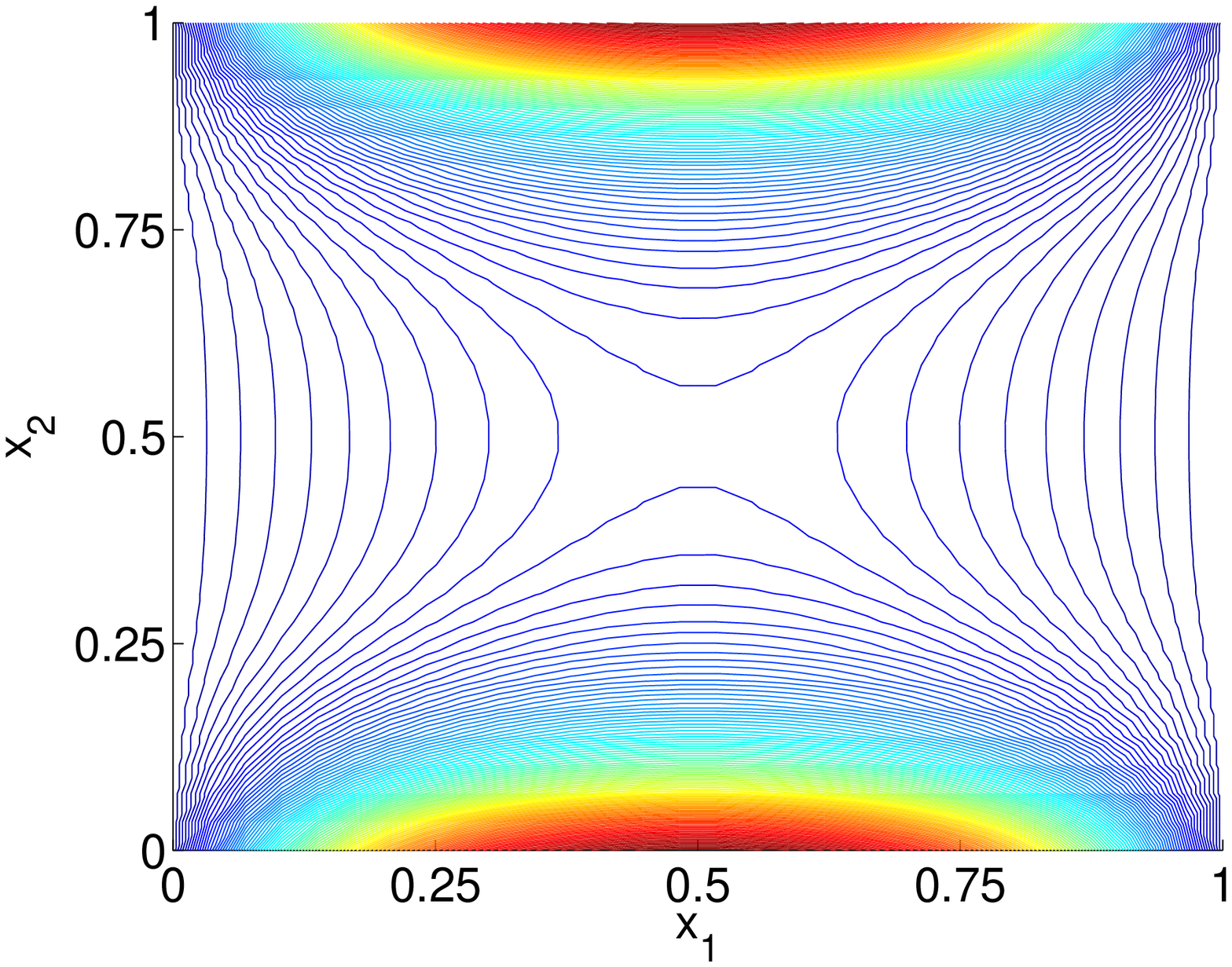}}
 \hspace{0.5cm}
 {\includegraphics[width=6.cm,height=6cm]{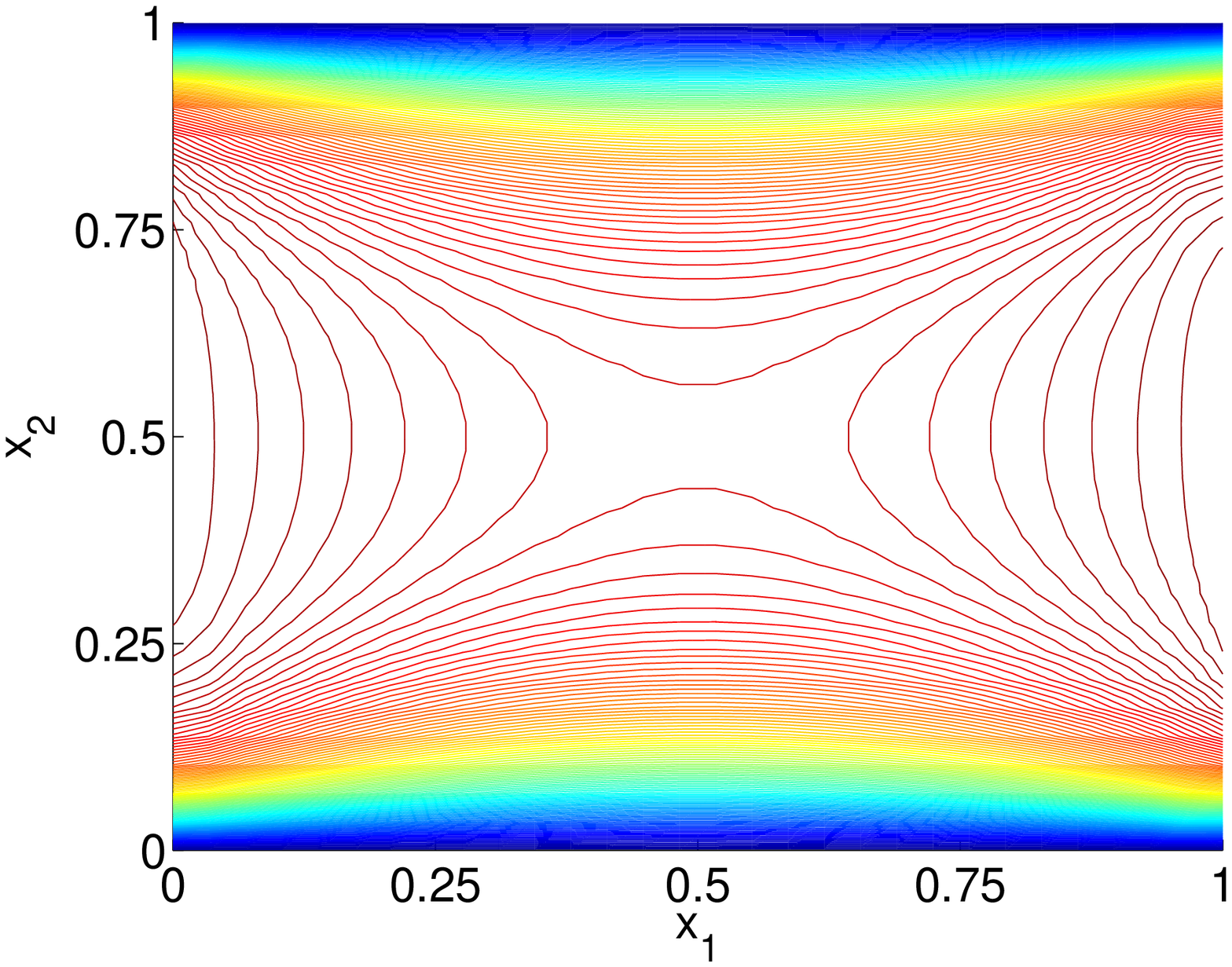}}
 \\ 
 {\includegraphics[width=8.cm,height=6.cm]{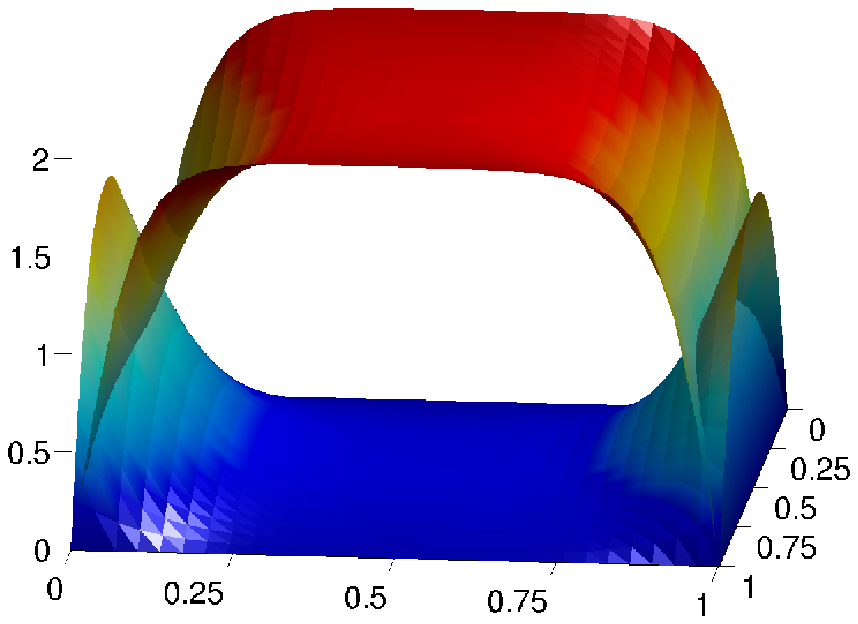}}
 %
 \caption{
Experiment 1. Top left is the level plot for $u_1$ and top right
is the level plot for $u_2$.
The bottom and top surfaces show $u_1$ and $u_2$ respectively.}
\label{fig.exp1}
\end{figure}

Recall that the diffusion coefficients, $c_i$ are
associated with growth rate of $u_i$.
In this experiment we implement (\ref{weak}) for different diffusion coefficients 
$c_1=0.1$, $c_2=0.01$ in (\ref{weak}) and $\eps=0$ in (\ref{weakeps}).
See Table  \ref{table.parameter} for parameter value.
Figure \ref{fig.exp1}
reflect the outcome of the distributions $u_1$ competing with
$u_2$.
%
%

Here $u_1$ shows effect similar to that shown for a single distribution of $u$ in 
\cite{Cohen2013}.
On the one hand, the distribution $u_1$ with $c_1 = 0.1$
ends up with higher density close to the boundaries compared to
the center of the values of its adaptive trait, $x_1$.
On the other, $u_2$ reveals trends over their values of adaptive traits
$x_2$ that mirror those of $u_1$.
These interactions among distributions of competing populations are reminiscent
of the principle of competitive exclusions,
where the population with slow diffusion---which is associated with 
slow growth rate---''winds''.

\subsubsection{Experiment 2: $c_1 = 0.1, c_2 = 0.1, \varepsilon = 0$}

\begin{figure}[t]
\centering
{\includegraphics[width=6.cm,height=6cm]{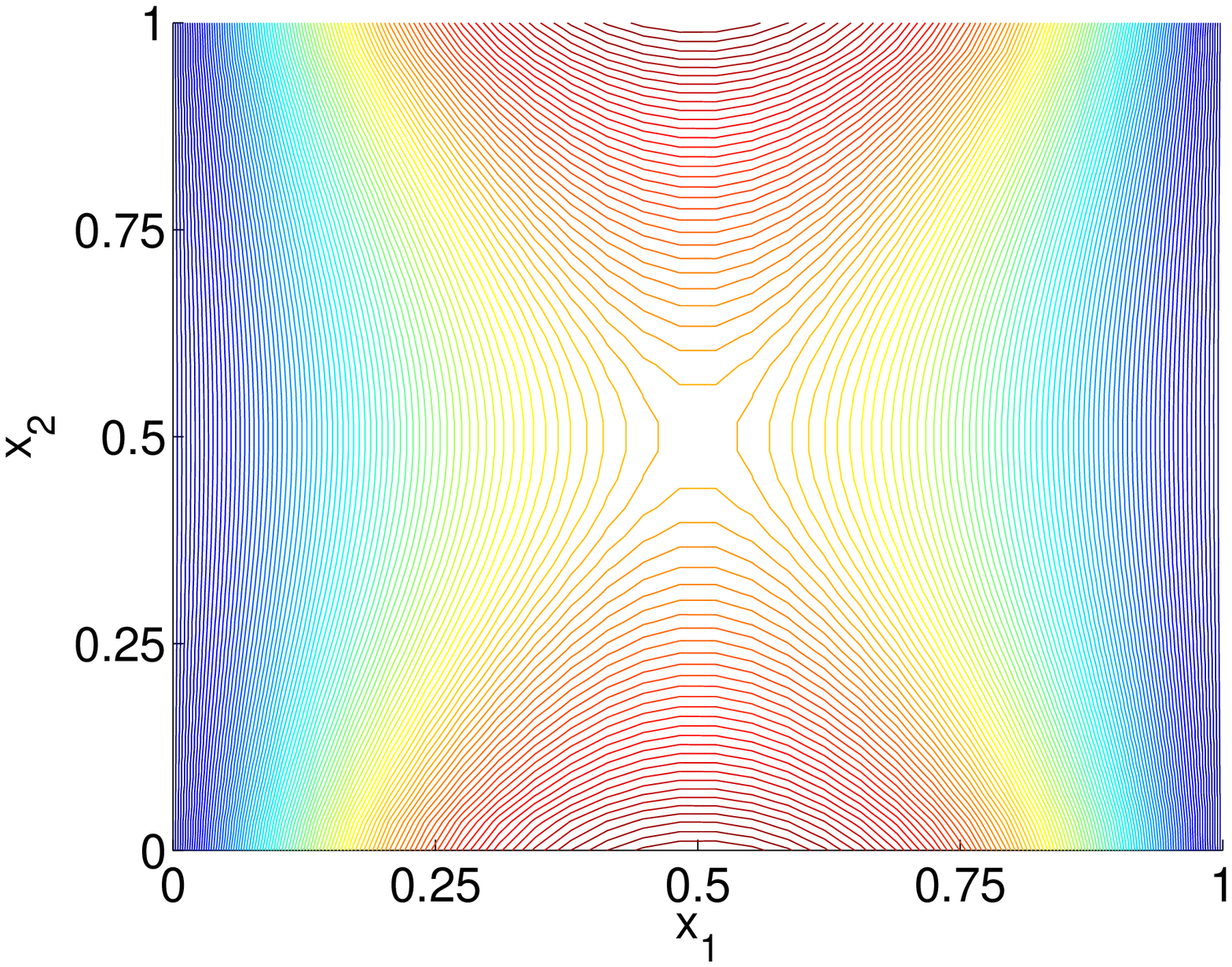}}
 \hspace{0.5cm}
 {\includegraphics[width=6.cm,height=6cm]{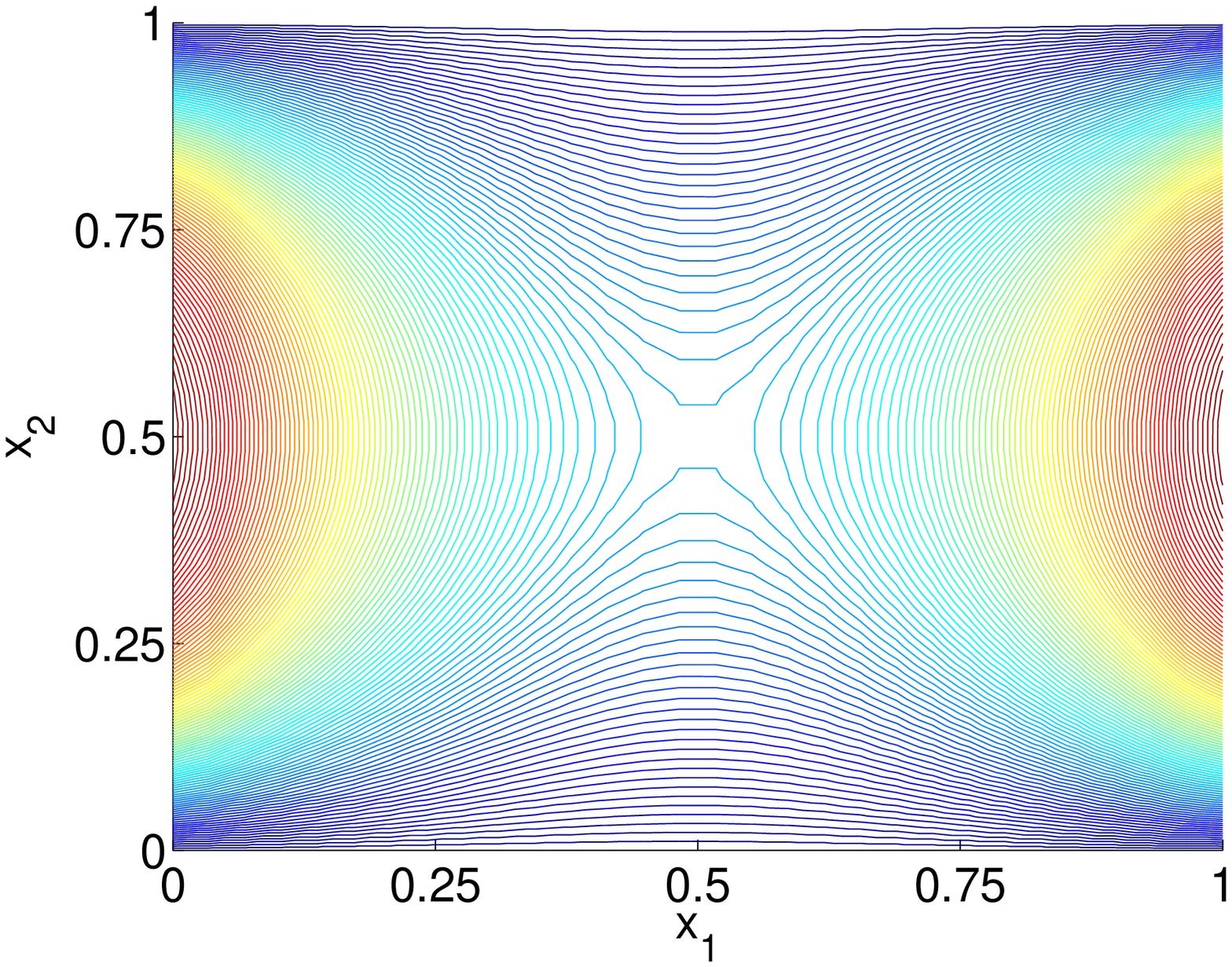}}
  {\includegraphics[width=8cm,height=6cm]{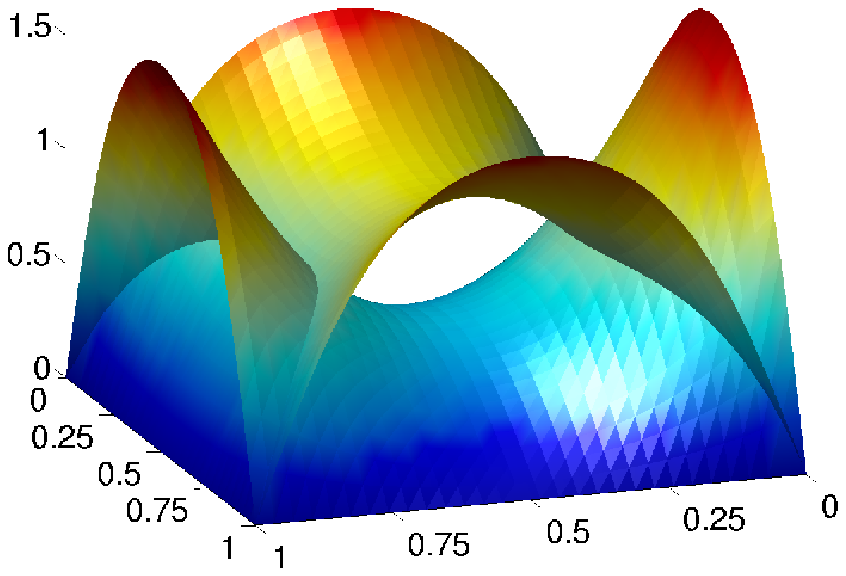}}
 \caption{Experiment 2. Top left is the level plot for $u_1$ and top right
is the level plot for $u_2$.
The bottom and top surfaces show $u_1$ and $u_2$ respectively.} 
\label{fig.exp2}
\end{figure}

Here we use $c_1=0.1$, $c_2=0.1$, $\eps=0$ (See Fig. \ref{fig.exp2}).
As expected, the distributions $u_1$ and $u_2$ mirror each other
in their density over the values of their adaptive traits $x_1$ and $x_2$.
Such symmetric distributions of phenotypes in two
populations is not likely to exist in nature.

\subsubsection{Experiment 3:  $c_1 = 0.1, c_2 = 0.1, \varepsilon > 0$ and
various boundary conditions}

The aim of this experiment is to investigate numerical results
for the singular perturbation problem \fer{weakeps}. 

The proof of Theorem~\ref{th.existence} implies the 
convergence of the sequence of solutions $(u_1^{(\eps)},u_2^{(\eps)})$ of \fer{weakeps}
to a solution of Problem P1.
The boundary conditions of problem \fer{weakeps} were set to the Dirichlet conditions.

\begin{figure}[t]
\centering
 {\includegraphics[width=6.25cm,height=4.75cm]{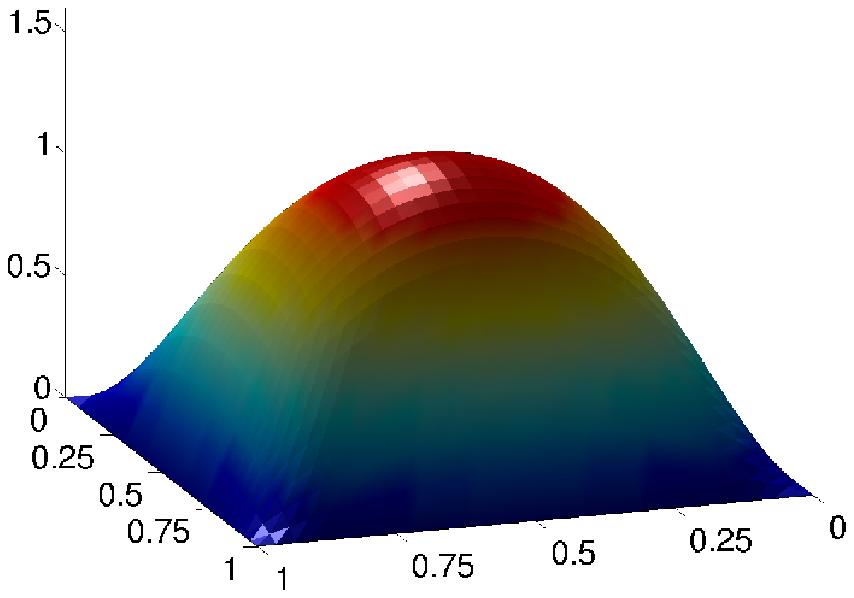}}
 \hspace{0.5cm}
 {\includegraphics[width=6.25cm,height=4.75cm]{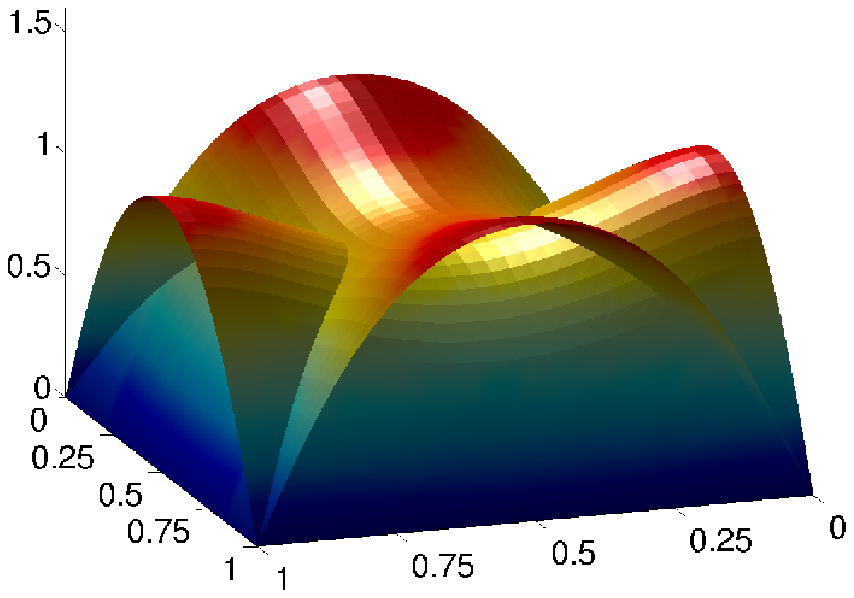}}
 \caption{{ Experiment 3. $\eps=0.1$. Left: Dirichlet. Right: Mixed}} 
\label{fig.exp3a}
\end{figure}

\begin{figure}[t]
\centering
 {\includegraphics[width=6.25cm,height=4.75cm]{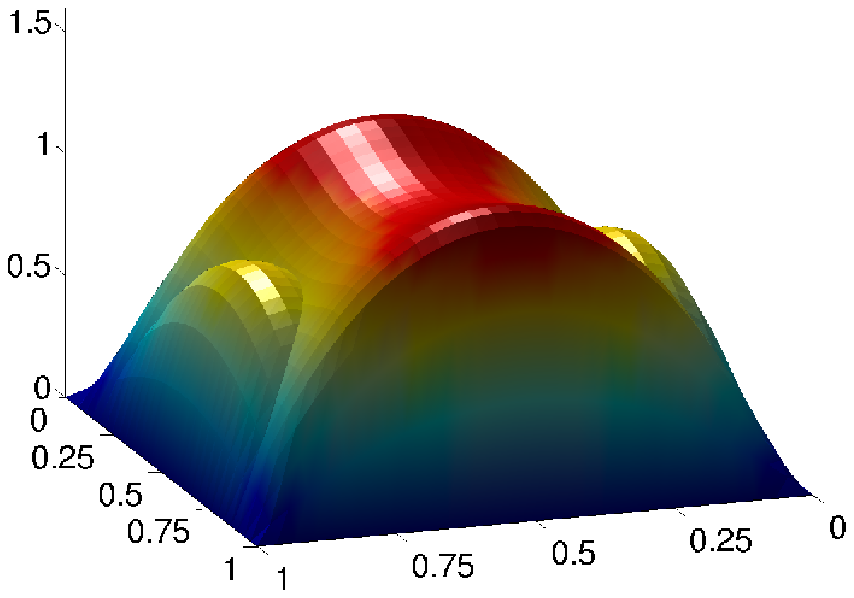}}
 \hspace{0.5cm}
 {\includegraphics[width=6.25cm,height=4.75cm]{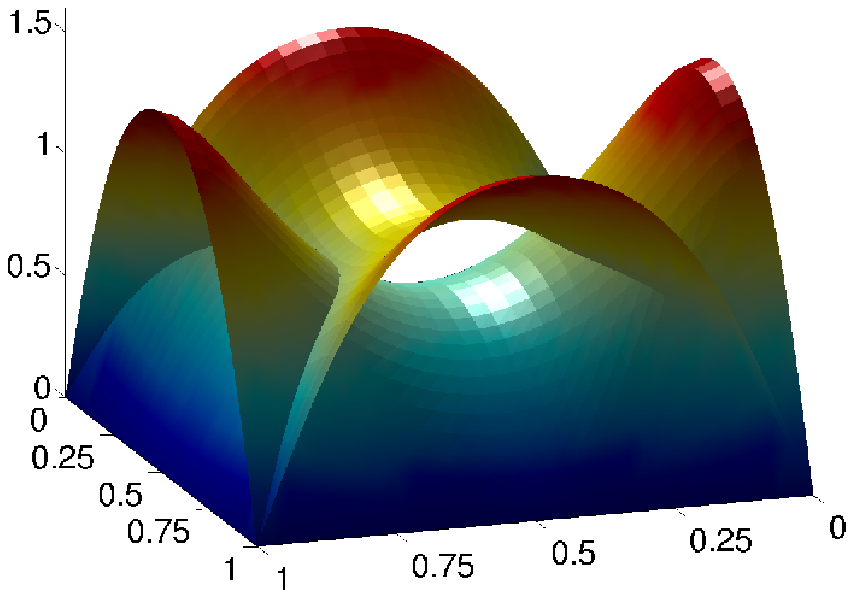}}
 \caption{{ Experiment 3. $\eps=0.01$. Left: Dirichlet. Right: Mixed}} 
\label{fig.exp3b}
\end{figure}

\begin{figure}[t]
\centering
 {\includegraphics[width=6.25cm,height=4.75cm]{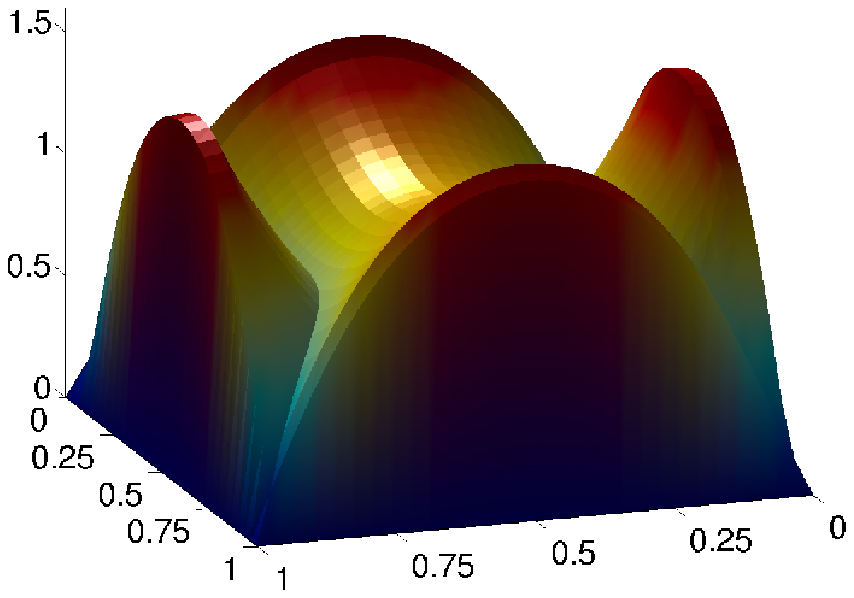}}
 \hspace{0.5cm}
 {\includegraphics[width=6.25cm,height=4.75cm]{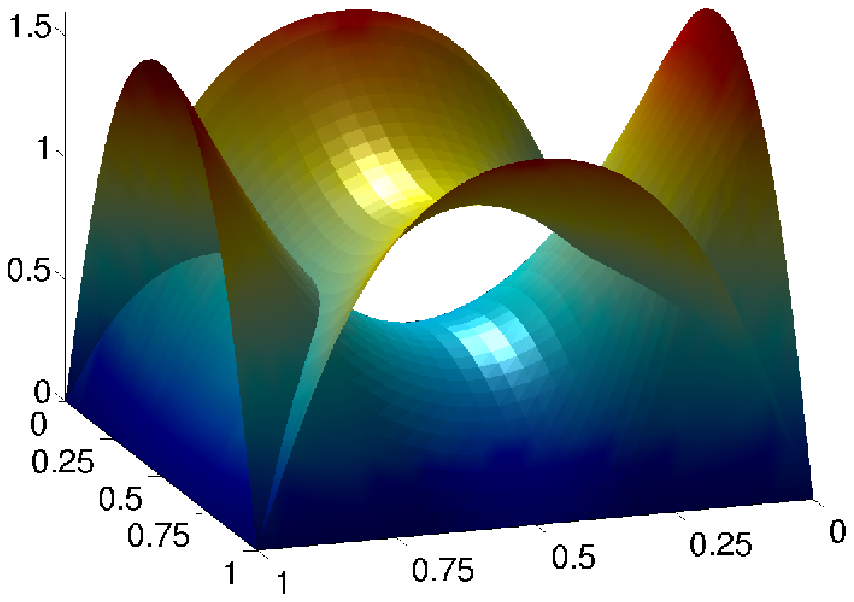}}
 \caption{{ Experiment 3. $\eps=1.e-10$. Left: Dirichlet. Right: Mixed}} 
\label{fig.exp3c}
\end{figure}

However, as  Figs. \ref{fig.exp3a}, \ref{fig.exp3b} and \ref{fig.exp3c} demonstrate, 
the mixed boundary conditions approximate the actual behavior of the 
unperturbed problem P1. By mixed boundary conditions we mean
\begin{align*}
 u_i^{(\eps)}=0 & \quad\text{on }(0,T)\times\partial \O_i\times \O_j\\
 \partial_{x_j} u_i^{(\eps)} =0 & \quad\text{on }(0,T)\times \O_i\times \partial\O_j
\end{align*}
for $i,j=1,2$, $i\neq j$.

In Figs. \ref{fig.exp3a}, \ref{fig.exp3b} and \ref{fig.exp3c} we show the steady state 
solutions for several values of $\eps$ for the singular perturbation problem 
with Dirichlet boundary conditions -DBC- (left panels) and mixed boundary conditions -MBC- (right panels). 
These solutions approximate those of Experiment~2 (see Fig.~\ref{fig.exp2}).

For $\eps=0.1$, the MBC approximate the unperturbed problem while the DBC force the 
solution to flatten near the boundary of $\O$, see Fig. \ref{fig.exp3a}.

For $\eps=0.01$, the MBC approximation still shows a declining in the mixed boundary, due to 
the Neumann boundary condition imposed in these boundaries. The DBC approximation is far from 
resembling the unperturbed solution (see Fig. \ref{fig.exp3b}).

Finally, for $\eps=10^{-10}$, the MBC approximation is indistinguishable from the solution to the unperturbed problem. 
However, the DBC approximation exhibits a steep boundary layer close to the mixed boundaries which only may 
disappear in the limit $\eps\to 0$, (see Fig. \ref{fig.exp3c}).

\section{Discussion}

Of all the disciplines in the biological sciences,
Evolution by Natural Selection is perhaps the
only one that can provide answers as to the question
of ``why'' as opposed to ``how''.

The dynamics of evolution provide unique challenges to its mathematical
modeling through PDE---which when applied to coevolving
populations we call a system of ED--in that (among others)
(\emph{i}) diffusion terms in the resulting PDEs are anisotropic and
(\emph{ii}) are determined, solely,
by the functions that model growth of populations.
Anisotropy in a system of PDEs refers to the property where
the diffusion term in each PDE is determined by said PDE only.

As opposed to single ED, 
the anisotropic property of ED surfaces only in systems of ED.
We  prove  existence of solutions such systems.

Because the structure of systems of ED is unique,
numerical algorithms to solve such systems need to be developed;
and we did.

Numerical experiments with such algorithms behave as expected from the latter.
We conduct factorial experiments that illustrate the effects of the perturbation parameter
$\varepsilon$ $\emph{vis \'a vis}$ Dirichlet and mixed boundary conditions.

The Dirichlet boundary conditions apply to models of ED.
This is so because phenotypes that carry values of adaptive traits
outside the boundaries of such traits cannot survive or cannot exist.
For example---except for few exceptional cases---organisms whose 
body temperature is regulated cannot
survive if it is (approximately) below $35^\circ$C.
Neither can they survive for long when their body temperature is (approximately)
above $42^\circ$C.

The numerical experiments with two competing populations,
each with its own adaptive trait, reveal the following:
(\emph{i}) The density of phenotypes mirror each other---this, 
in spite of the anisotropy; 
such a pattern surfaces because the ED of the populations
interact through the competition term.
(\emph{ii}) Near the boundaries, and as expected,
the population that diffuses faster ($c_1 = 0.1$)
obtains density higher than the population that diffuses
slower ($c_1 = 0.01$).
For a single ED, we \cite[]{Cohen2013} obtained  dynamics that
correspond to the ED with $c_1 = 0.1$ (see left and bottom panels
in Figure \ref{fig.exp1}) regardless of the magnitude of the diffusion coefficient.
This is so because in a system with a single ED there is no influence
(in terms of the diffusion) of
another ED that prevents it from achieving high density near the boundaries.

Because $\varepsilon = 0$, the results of experiment $2$ reveal symmetric
distributions of $u_1$ and $u_2$; each of the ED correspond to our finding
for a single ED \cite[]{Cohen2013}.

The results from experiments $3$ to $5$
(given in Figures \ref{fig.exp3a} - \ref{fig.exp3c}) 
illustrate the influence of the perturbation parameter
$\varepsilon$ with Dirichlet boundary conditions compared to mix boundary
conditions.
The right panels of each of the Figures reflect the ever increasing density
at the boundaries (compared to the center) of the distribution.
The trends in these results correspond to those with the Dirichlet
boundary conditions.


\end{document}